\documentclass{amsart}
\begin{document}
\title{Note on the question of Sikora}
\author{K. V. Storozhuk}
\address{Sobolev Institute of Mathematics, Acad.
Koptyug pr. 4, 630090 Novosibirsk, Russia}


\begin{abstract} A natural topology on the set of left orderings on free
   abelian groups and free groups $F_n$, $n>1$ has studied in [1].
   It has been proven already that in the abelian case the resulted
    topological space is a Cantor set. There was a conjecture:
    this is also true for the free group $F_n$ with $n>1$ generators.
    We point out the article dealing with equivalent questions.

\end{abstract} \maketitle

Following [1], denote the set of left (right) orderings on a group
 $G$ by $LO(G)$ $(RO(G)$); put $BiO(G)=LO(G)\cap RO(G)$. In [1], a
 topology is defined on the set of orderings. The base of this
 topology is the collections of orderings defined by a finite number
 of inequalities.

 This topology is conceptually close to the ideal topology of a ring.

 For a countable group $G$, the space of orderings is a compact,
 totally disconnected metric space. Thus, for this space to be
 homeomorphic to a Cantor set it suffices that the following
 condition holds ([1],Corollary 1.6): each set from the base of
 topology is either empty or infinite.

 It is proven in [1] that the sets $LO(G)$ and $BiO(G)$ are Cantor
 sets for $G=Z^n$, and a conjecture is put forward that this sets for
 $G=F_n$ are Cantor sets as well.
 Note that this question (for left orders) and some corollaries answered in
[4]

 Notice that, in our context, the discussed requirement (corollary
 1.6) is equivalent to the absence of isolated points. Now, the
 conjecture of [1] can be algebraically reformulated as follows: does
 there exist a right ordering (two-sided ordering) on a free group
 $F_n$ that is defined uniquely by a finite collection of
 inequalities?

 It is in this fashion that the question is posed in [2]. For
 right orderings, the negative answer is announced there and proven
 in [3]. So, the discussed conjecture is true for $RO(F_n)$ (and, of
 course, for $LO(F_n)$).
  The question for $BiO(F_n)$ seems to be
 still open --- together with the corresponding question from
 [2].

 \vskip1mm

The author is grateful to Professor V. M. Kopytov, who drew the
 author's attention to the article [3].

   \email{stork@math.nsc.ru}

\end{document}